\documentclass[12pt]{article}
\usepackage{amsmath, amssymb}
\begin{document}
\centerline{\bf Finite Generation of Canonical Ring by Analytic
Method}

\bigbreak\centerline{\it Dedicated to Professor Lu Qikeng on his
80th Birthday} \bigbreak
\centerline{Yum-Tong Siu\ %
\footnote{Partially supported by a grant from the National Science
Foundation.} }

\bigbreak\noindent{\bf \S0.} {\sc Introduction.}  The analytic
methods of $L^2$ estimates of $\bar\partial$ and multiplier ideal
sheaves provide a powerful new approach to a number of long
outstanding problems in algebraic geometry.  Besides effective
results on problems related to the Fujita conjecture and the
Matsusaka big theorem, the deformational invariance of the
plurigenera was proved by such an approach first for the case of
general type [Siu 1998] and then for the general algebraic case [Siu
2002].

\medbreak The techniques developed for the deformational invariance
of the plurigenera were intended to prove the finite generation of
the canonical ring.  The extension result on pluricanonical sections
from the method of the deformational invariance of the plurigenera
opens up the possibility of using restriction to hypersurfaces and
induction on dimension to prove the finite generation of the
canonical ring.  However, from an analytic viewpoint the technical
details arising from the various singular situations are quite
daunting.

\medbreak The deformational invariance of the plurigenera was proved
by using the techniques of (i) the global generation of multiplier
ideal sheaves (A.1), (ii) the extension theorem of Ohsawa-Takegoshi
[Ohsawa-Takegoshi 1987] (which in this setting can be replaced by
the vanishing theorem for multiplier ideal sheaves [Kawamata 1982,
Viehweg 1982, Nadel 1990] ), and (iii) the use of one of the
canonical bundle inside the pluricanonical bundle as the volume form
to be used in the $L^2$ estimates of $\bar\partial$. The proof of
the finite generation of the canonical ring can more easily be
handled by directly applying the above three techniques from the
proof of the deformational invariance of the plurigenera.

\medbreak In this article we give an overview of the analytic proof
of the following theorem on the finite generation of the canonical
ring for the case of general type.

\bigbreak\noindent(0.1) {\it Theorem.}  Let $X$ be a compact complex
algebraic manifold of general type.  Then the canonical ring
$$
R\left(X, K_X\right)=\bigoplus_{m=1}^\infty\Gamma\left(X,
mK_X\right)
$$
is finitely generated.

\bigbreak\noindent Details of techniques for the analytic proof of
the finite generation of the canonical ring for the case of general
type were posted in [Siu 2006, Siu 2007].  An algebraic proof was
posted in [Birkan-Cascini-Hacon-McKernan 2006].

\medbreak In this overview we focus on the formulation using the
notion of a discrepancy subspace, which measures the extent of
failure of achieving stable vanishing order in terms of uniformity
in $m$ for all $m$-canonical bundles.  It highlights more clearly
how the analytic method handles the problem of infinite number of
interminable blow-ups in the intuitive approach to prove the finite
generation of the canonical ring.

\medbreak Toward the end of this overview we discuss how our
situation is similar to what is needed for a proof of the abundance
conjecture.  An adaptation of the argument here for use in a proof
of the abundance conjecture would require an analytic argument of
controlling the estimates in passing to limit, which is analogous to
the situation of extending the proof of the deformational invariance
of the plurigenera for the case of general type [Siu 1998] to the
general algebraic case without the general type assumption [Siu
2002].

\medbreak The notations ${\mathbb C}$, ${\mathbb Q}$ and ${\mathbb
Z}$, and ${\mathbb N}$ denote respectively the complex numbers, the
rational numbers, the integers, and the positive integers.  The
reduced structure sheaf of a complex space $W$ is denoted by
${\mathcal O}_W$.  The maximum ideal of a point $P$ of a complex
manifold $Y$ is denoted by ${\mathfrak m}_P$.  The ideal sheaf of a
subvariety $Z$ in a complex manifold $Y$ is denoted by ${\mathcal
I}_Z$.  The canonical section of the line bundle associated to a
complex hypersurface $V$ in a complex manifold $Y$ is denoted by
$s_V$.  The multiplier ideal sheaf of a metric $e^{-\psi}$ is
denoted by ${\mathcal I}_\psi$ (see (A.7)). The canonical line
bundle of a complex manifold $Y$ is denoted by $K_Y$. A multi-valued
holomorphic section $\sigma$ of a holomorphic line bundle $E$ over a
complex manifold $Y$ means that $\sigma^m$ is a holomorphic section
of $mE$ over $Y$ for some positive integer $m$.

\bigbreak\noindent{\bf\S 1.} {\sc Setup of Descrepancy Subspace.} In
this overview of the analytic proof of the finite generation of the
canonical ring, for the sake of clarity and transparency of the core
arguments we skip all arguments involving diophantine approximations
so that certain real numbers which need to be proved by diophantine
arguments to be rational are just assumed to be known to be
rational.  In order not to disrupt the main line of the core
arguments, we move to the Appendix many side arguments and the
listing of known facts and techniques together with some simple
adaptations needed for our purpose.  Sometimes the condition of a
positive integer being sufficiently divisible is simply stated
sloppily as being sufficiently large.

\bigbreak\noindent(1.1) {\it Metric of Minimum Singularity and Its
Truncation.}  Let $X$ be a compact complex algebraic manifold of
general type.  Let
$$
s^{(m)}_1, \cdots, s^{(m)}_{q_m}\in\Gamma\left(X, mK_X\right)
$$
be a ${\mathbb C}$-linear basis and let
$\left\{\varepsilon_m\right\}_{m\in{\mathbb N}}$ be a sequence of
positive constants decreasing so rapidly monotonically to zero that
$$
\Phi=\sum_{m=1}^\infty\varepsilon_m\sum_{j=1}^{q_m}\left|s^{(m)}_j\right|^{\frac{2}{m}}
$$
converges on $X$.  Let $\varphi=\log\Phi$ so that
$\frac{1}{\Phi}=e^{-\varphi}$ is a metric for $K_X$ which we call
{\it the metric of $K_X$ of minimum singularity}.  (It actually is
not unique, but depends on the choice of the basis $s^{(m)}_1,
\cdots, s^{(m)}_{q_m}$ and the sequence
$\left\{\varepsilon_m\right\}_{m\in{\mathbb N}}$ of positive
numbers.)  For $N\in{\mathbb N}$ let
$$
\Phi_N=\sum_{m=1}^N\varepsilon_m\sum_{j=1}^{q_m}\left|s^{(m)}_j\right|^{\frac{2}{m}}
$$
We call the metric $\frac{1}{\Phi_N}$ for $K_X$ the {\it $N$-th
truncation of the metric of minimum singularity}.

\medbreak\noindent(1.2) {\it Stable Vanishing Order.}  For a point
$P$ of $X$ we say that the stable vanishing order is achieved at $P$
if there exists some $N\in{\mathbb N}$ such that the two functions
$\Phi_N$ and $\Phi$ are comparable on some open neighborhood $U$ of
$P$ in $X$ in the sense that $\Phi\leq C\Phi_N$ for some positive
constant $C$ on $U$.

\medbreak When we say that the generic stable vanishing order across
some subvariety $Y$ of $X$ is achieved, we mean that at a generic
point of $Y$ the vanishing order of $\Phi$ across $Y$ is the same as
the vanishing order of $\Phi_N$ across $Y$ for some $N\in{\mathbb
N}$.

\medbreak It follows from a direct application of Skoda's result on
ideal generation that when the stable vanishing order is achieved at
every point of $X$, the canonical ring is finitely generated (see
Appendix (A.4)). So the problem of the finite generation of the
canonical ring is reduced to proving the achievement of stable
vanishing order everywhere on $X$.

\medbreak The most natural approach is to consider the subvariety
$Z$ of $X$ where the stable vanishing order is not achieved and then
try to prove that if $Z$ is nonempty, we can show that at one of its
points the stable vanishing order is achieved, giving a
contradiction.

\medbreak However, it turns out that in deriving a contradiction for
not achieving the stable vanishing order everywhere, it is more
efficient to get a contradiction for a statement incorporating the
extent (or multiplicity) of the failure of the achievement of the
stable vanishing order.  We introduce now the notion of discrepancy
subspace which measures the extent (or multiplicity) of the failure
of the achievement of the stable vanishing order.

\bigbreak\noindent(1.3) {\it Definition of Discrepancy Subspace.}
Let ${\mathcal J}$ be a coherent ideal sheaf on $X$. The stable
vanishing order of the canonical line bundle of $X$ is said to be
{\it precisely achieved modulo the subspace of $X$ defined by
${\mathcal J}$} if there exist some positive integer $m_{\mathcal
J}$ and some positive constant $C_{m,k,{\mathcal J}}$ for $k,
m\in{\mathbb N}$ with $m\geq m_{\mathcal J}$ such that the
inequality
$$
\left|s_{\mathcal
J}\right|^2\,\sum_{j=1}^{q_k}\left|s^{(k)}_j\right|^{\frac{2m}{k}}\leq
C_{m,k,{\mathcal
J}}\sum_{j=1}^{q_m}\left|s^{(m)}_j\right|^2\leqno{(\dagger)}
$$
holds on $X$ for all $k, m\in{\mathbb N}$ with $m\geq m_{\mathcal
J}$. Here the notation $\left|s_{\mathcal J}\right|^2$ means the
following. For a coherent ideal sheaf ${\mathcal I}$ on $X$
generated locally by holomorphic function germs
$\tau_1,\cdots,\tau_\ell$, we define
$$
\left|s_{\mathcal I}\right|^2=\sum_{j=1}^\ell\left|\tau_j\right|^2.
$$
Let $Z$ be the zero-set of the coherent ideal sheaf ${\mathcal J}$.
We call the ringed space $\left(Z,{\mathcal O}_X\left/{\mathcal
J}\right.\right)$ a {\it discrepancy subspace}.  We call the
coherent ideal sheaf ${\mathcal J}$ a {\it discrepancy ideal sheaf}.

\bigbreak\noindent(1.3.1) {\it Remark.} The main feature of the
formulation is that a discrepancy subspace not only specifies the
set of points (which is the zero-set of the ideal sheaf ${\mathcal
J}$) where the failure of achieving the stable vanishing order
occurs but also describes the extent of the failure by providing the
ideal sheaf ${\mathcal J}$, simultaneously for all $mK_X$ for
sufficiently large $m$, so that if one adds the vanishing order of
this ideal sheaf the stable vanishing order is no less than that
given by a finite number of pluricanonical sections. The vanishing
order of an ideal sheaf is formulated here in terms of the sum of
the absolute-value-squares of its local generators.  This
description of the extent of the failure by the ideal sheaf
${\mathcal J}$ makes the induction process easier to handle.  The
simultaneity for all $mK_X$ for sufficiently large $m$ holds the key
to understanding the reason for the termination of the infinite
process of blow-ups in the intuitive approach of blowing up
successively to prove the finite generation of the canonical ring.

\bigbreak\noindent(1.4) {\it Formulation in Terms of Metric of
Minimum Singularity.} In the formulation $(\dagger)$ the constant
$C_{m,k,{\mathcal J}}$ is allowed to depend on $m$ and $k$, because
we are free to choose the rapidly decreasing sequence of positive
numbers $\varepsilon_\ell$ used in the definition of
$$
\Phi=\sum_{\ell=1}^\infty\varepsilon_\ell\sum_{j=1}^{q_\ell}\left|s^{(\ell)}_j\right|^{\frac{2}{\ell}}
$$
in order to obtain the inequality $\Phi\leq C\Phi_N$ for some
positive number $C$ and some positive integer $N$ with
$$
\Phi_N=\sum_{\ell=1}^N\varepsilon_\ell\sum_{j=1}^{q_\ell}\left|s^{(\ell)}_j\right|^{\frac{2}{\ell}}.
$$
Another way to formulate $(\dagger)$ is that for any fixed $m\geq
m_{\mathcal J}$ and for an $m$-dependent appropriate choice of the
rapidly decreasing sequence of positive numbers $\varepsilon_\ell$
used in the definition of
$$
\Phi=\sum_{\ell=1}^\infty\varepsilon_\ell\sum_{j=1}^{q_\ell}\left|s^{(\ell)}_j\right|^{\frac{2}{\ell}}
$$
the inequality
$$
\left|s_{\mathcal J}\right|^2\,\Phi^m\leq\tilde C_{m,{\mathcal
J}}\sum_{j=1}^{q_m}\left|s^{(m)}_j\right|^2\leqno{(\dagger)_m^\flat}
$$
holds on $X$ for some constant $\tilde C_{m,{\mathcal J}}$.  In
other words, we use the inequality
$$
\left|s_{\mathcal J}\right|^2\left(\check{\Phi}_m\right)^m\leq\tilde
C_{m,{\mathcal J}}\sum_{j=1}^{q_m}\left|s^{(m)}_j\right|^2,
$$
where
$$
\check{\Phi}_m=\sum_{k=1}^\infty\varepsilon_{k,m}\sum_{j=1}^{q_k}\left|s^{(k)}_j\right|^{\frac{2}{k}}
$$
for some positive constants $\varepsilon_{k,m}$.

\medbreak\noindent(1.5) {\it Transformation of Discrepancy Subspace
in Blow-Up and Conductor.}  In a blow-up $\tilde X\to X$ of $X$, the
r\^ole played by the adjunction formula is canceled by its effect on
both sides of $(\dagger)$ and $\left|s_{\mathcal J}\right|^2$ simply
transforms as a lifting of a local function from $X$ to $\tilde X$.
This enables us to assume that, after replacing $X$ by its blowup,
${\mathcal J}$ is the ideal sheaf of a divisor whose components are
in normal crossing.  With the blow-up, we can use the technique of
the minimum center of log canonical singularity [Kawamata 1985,
Shokurov 1985].

\medbreak Also note that from the inequality
$$
\left|s_{\mathcal J}\right|^2\,\Phi^m\leq\tilde C_{m,{\mathcal
J}}\sum_{j=1}^{q_m}\left|s^{(m)}_j\right|^2
$$
in $(\dagger)_m^\flat$, if
$$
\int\frac{\left|f\right|^2}{\Phi^m}
$$
is finite for some local holomorphic function germ $f$, then
$$
\int\frac{\left|s_{\mathcal
J}\right|^2\left|f\right|^2}{\sum_{j=1}^m\left|s^{(m)}_j\right|^2}\leq
\tilde C_{m,{\mathcal J}}\int\frac{\left|f\right|^2}{\Phi^m}
$$
is also finite.  This means that
$$
{\mathcal J}{\mathcal I}_{m\log\Phi}\subset{\mathcal
I}_{\log\sum_{j=1}^m\left|s^{(m)}_j\right|^2}
$$
and the key point is that the {\it conductor} ${\mathcal J}$ is
independent of $m$.  In particular, we are able to locate the
minimum center of log canonical singularity (or irreducible subspace
of minimum discrepancy) in a way which works for all $mK_X$ with $m$
sufficiently large.

\bigbreak\noindent(1.6) {\it Intersection of Discrepancy Subspaces
and Minimum Discrepancy Subspace.} Though for the precise
achievement of stable vanishing order we need only show that
$\Phi^m$ and $\sum_{j=1}^{q_m}\left|s^{(m)}_j\right|^2$ are
comparable for some $m\in{\mathbb N}$, the reason why we need
$(\dagger)_m^\flat$ for all $m\geq m_{\mathcal J}$ and not just for
some single $m$ is that we use induction to reduce the subspace
defined by ${\mathcal J}$ step-by-step, which means that by
replacing $m_{\mathcal J}$ by an appropriately larger
$m_{\tilde{\mathcal J}}$ in the comparison between $\Phi^m$ and
$\sum_{j=1}^{q_m}\left|s^{(m)}_j\right|^2$ we seek to replace
${\mathcal J}$ by a strictly bigger ideal sheaf $\tilde{\mathcal J}$
in the inequality $(\dagger)_m^\flat$.  Since we need to use a
bigger $m_{\mathcal J}$ in every one of the finite-step inductive
process, the inequality $(\dagger)_m^\flat$ has to be formulated to
hold for all $m\geq m_{\mathcal J}$. Another way to look at this is
as follows.  Suppose we have, not only
$$
\left|s_{\mathcal
J}\right|^2\,\sum_{j=1}^{q_k}\left|s^{(k)}_j\right|^{\frac{2m}{k}}\leq
C_{m,k,{\mathcal
J}}\sum_{j=1}^{q_m}\left|s^{(m)}_j\right|^2\leqno{(\dagger)_{\mathcal
J}}
$$
on $X$ for all $k, m\in{\mathbb N}$ with $m\geq m_{\mathcal J}$, but
also
$$
\left|s_{\tilde{\mathcal
J}}\right|^2\,\sum_{j=1}^{q_k}\left|s^{(k)}_j\right|^{\frac{2m}{k}}\leq
C_{m,k,\tilde{\mathcal
J}}\sum_{j=1}^{q_m}\left|s^{(m)}_j\right|^2\leqno{(\dagger)_{\tilde{\mathcal
J}}}
$$
on $X$ for all $k, m\in{\mathbb N}$ with $m\geq m_{\tilde{\mathcal
J}}$.  Then we can define ${\mathcal K}$ as the sum of ${\mathcal
J}$ and $\tilde{\mathcal J}$ and set $m_{\mathcal K}$ as the maximum
of $m_{\mathcal J}$ and $m_{\tilde{\mathcal J}}$ and set
$C_{m,k,{\mathcal K}}$ as the sum of $C_{m,k,{\mathcal J}}$ and
$C_{m,k,{\tilde{\mathcal J}}}$ and get
$$
\left|s_{\mathcal
K}\right|^2\,\sum_{j=1}^{q_k}\left|s^{(k)}_j\right|^{\frac{2m}{k}}\leq
C_{m,k,{\mathcal
K}}\sum_{j=1}^{q_m}\left|s^{(m)}_j\right|^2\leqno{(\dagger)_{\mathcal
K}}
$$
on $X$ for all $k, m\in{\mathbb N}$ with $m\geq m_{\mathcal K}$.

\medbreak On the other hand, if we have only
$$
\left|s_{\mathcal J}\right|^2\,\Phi^m\leq\tilde C_{m,{\mathcal
J}}\sum_{j=1}^{q_m}\left|s^{(m)}_j\right|^2
$$
for some constant $\tilde C_{m,{\mathcal J}}$ for one single
$m=m_{\mathcal J}$, then having $$\left|s_{\tilde{\mathcal
J}}\right|^2\,\Phi^m\leq\tilde C_{m,\tilde{\mathcal
J}}\sum_{j=1}^{q_m}\left|s^{(m)}_j\right|^2
$$
for some constant $\tilde C_{m,\tilde{\mathcal J}}$ for one single
$m=m_{\tilde{\mathcal J}}$ would not be sufficient to give us
$$
\left|s_{\mathcal K}\right|^2\,\Phi^m\leq\tilde C_{m,{\mathcal
K}}\sum_{j=1}^{q_m}\left|s^{(m)}_j\right|^2
$$
for some constant $\tilde C_{m,{\mathcal K}}$ for one single
$m=m_{\mathcal K}$ when ${\mathcal K}$ is the sum of ${\mathcal J}$
and $\tilde{\mathcal J}$.  Even, if we take the least common
multiple $\check{m}$ of $m_{\mathcal J}$ and $m_{\tilde{\mathcal
J}}$ so that $\check{m}=m_{\mathcal J} p_{\mathcal J}$ and
$\check{m}=m_{\tilde{\mathcal J}} p_{\tilde{\mathcal J}}$, then we
can only get
$$
\left|s_{\check{\mathcal K}}\right|^2\,\Phi^{\check{m}}\leq\check
C\sum_{j=1}^{q_{\check{m}}}\left|s^{\left(\check{m}\right)}_j\right|^2
$$
for some $\check C>0$, where $\check{\mathcal K}$ is generated by
${\mathcal J}^{p_{\mathcal J}}$ and ${\tilde{\mathcal
J}}^{p_{\tilde{\mathcal J}}}$.  In this case, in general
$\check{\mathcal K}$ is not equal to ${\mathcal K}$.

\medbreak Since we can take the sum of two coherent ideal sheaves
for discrepancy subspaces and since we have the Noetherian property
for a nondecreasing chain of coherent ideal sheaves, we can get a
{\it minimum discrepancy subspace} defined by a {\it maximum
discrepancy ideal sheaf}. Moreover, by blowing up, we can assume
that the maximum discrepancy ideal sheaf is the ideal sheaf for a
divisor whose components are nonsingular hypersurfaces in normal
crossing.

\bigbreak\noindent(1.7) {\it Non-Achievement of Stable Vanishing
Order and Discrepancy Subspace.} Suppose ${\mathcal J}$ is the
maximum discrepancy ideal sheaf and ${\mathcal J}$ is equal to the
ideal sheaf of the divisor
$$
D=\sum_{j=1}^\ell\alpha_j D_j
$$
with $\left\{D_j\right\}_{1\leq j\leq\ell}$ composed of nonsingular
hypersurfaces in normal crossing and $\ell\in{\mathbb N}$ and
$\alpha_j\in{\mathbb N}$ for $1\leq j\leq\ell$.  It does not mean
that at points of $D_1-\cup_{j=2}^\ell D_j$ the generic stable
vanishing order definitely cannot be achieved across $D_1$.  The
reason is the following.  For simplicity let us explain this in the
case of $\ell=1$.  Logically it may happen that there are coherent
ideal sheaves ${\mathcal I}_m$ on $X$ for $m\in{\mathbb N}$ whose
zero-set $Z_m$ is a proper subvariety of $D_1$ and the inequality
$$
\left|s_{{\mathcal I}_m}\right|^2\left(\check{\Phi}_m\right)^m\leq
C^\sharp_m\sum_{j=1}^{q_m}\left|s^{(m)}_j\right|^2
$$
holds for some positive constant $C^\sharp_m$, where
$$
\check{\Phi}_m=\sum_{k=1}^\infty\varepsilon_{k,m}\sum_{j=1}^{q_k}\left|s^{(k)}_j\right|^{\frac{2}{k}}
$$
for some positive constants $\varepsilon_{k,m}$, but ${\mathcal
I}_{m+1}$ is a proper ideal subsheaf of ${\mathcal I}_m$ for
$m\in{\mathbb N}$ so that there is no coherent ideal ${\mathcal I}$
on $X$ independent of $m$ such that ${\mathcal I}$ is contained in
each ${\mathcal I}_m$ for $m\in{\mathbb N}$ and the zero-set of
${\mathcal I}$ is a proper subvariety $Z$ of $D_1$.  What may happen
is that ${\mathcal J}$ may be equal to the intersection
$\cap_{m=1}^\infty{\mathcal I}_m$ and, though the zero-set of each
${\mathcal I}_m$ is the proper subvariety $Z_m$ of $D_1$, yet the
zero-set of the intersection $\cap_{m=1}^\infty{\mathcal I}_m$ is
$D_1$.  In a way this formulation of using discrepancy subspaces
differs from the simple precise achievement of generic stable
vanishing order in that it removes the difficulty of the
ever-changing setting as the $m$ in $mK_X$ increases but it also
makes the task of decreasing the discrepancy subspace so much
harder.

\medbreak\noindent(1.7.1) {\it Remark.} When we try to decrease the
discrepancy ideal subspace, the situation may seem simpler if the
generic stable vanishing order is achieved at some point of $D$.
However, it is not exactly the case, because the discrepancy
measures the difference between bigness and some appropriate
ampleness and does not just measure the difference between bigness
and numerical effectiveness. When the generic stable vanishing order
$\gamma$ across $Y$ is achieved, the ${\mathbb Q}$-line bundle
$K_X-\gamma Y$ may be locally numerically effective at some affine
open subset $\Omega_Y$ of $Y$, yet for $K_X-\gamma Y$ at points of
$\Omega_Y$ there may still be a difference between the local
numerical effectiveness and the required ampleness at the points of
$\Omega_Y$.

\bigbreak\noindent(1.8) {\it Geometric Construction of Minimum
Discrepancy Subspace.} There is a more constructive way to identify
the minimum discrepancy subspace. For fixed $k$ and $m$ let
${\mathcal A}_{k,m}$ be the ideal sheaf consisting of all
holomorphic function germs $f$ on $X$ at a point $P$ of $X$ such
that
$$
\left|f\right|^2\,\sum_{j=1}^{q_k}\left|s^{(k)}_j\right|^{\frac{2m}{k}}\leq
C_{f,P}\sum_{j=1}^{q_m}\left|s^{(m)}_j\right|^2\leqno{(\dagger)_{k,m}}
$$
on some open neighborhood $U$ of $P$ for some positive constant
$C_{f,P}$.  By using blow-ups of $X$, it is clear that ${\mathcal
A}_{k,m}$ is coherent.  Let ${\mathcal B}_m$ be the largest coherent
ideal sheaf on $X$ such that ${\mathcal B}_m$ is contained in
$\cap_{k,\ell\geq m}{\mathcal A}_{k,\ell}$.  Clearly ${\mathcal
B}_m$ is contained in ${\mathcal B}_{m+1}$, because
$\cap_{k,\ell\geq m}{\mathcal A}_{k,\ell}$ is contained in
$\cap_{k,\ell\geq m+1}{\mathcal A}_{k,\ell}$.  It follows from the
Noetherian property of a convergent power series ring and the
compactness of $X$ that there exists some $m_0\in{\mathbb N}$ such
that ${\mathcal B}_{m_0}={\mathcal B}_m$ for $m\geq m_0$. The
coherent ideal sheaf ${\mathcal B}_{m_0}$ is equal to the maximum
discrepancy ideal sheaf ${\mathcal J}$ with $m_{\mathcal J}=m_0$.

\medbreak Note that the discrepancy subspace defined by ${\mathcal
J}$ must be inside $D$ when $aK_X=A+D$ (for some $a\in{\mathbb N}$
and some ample $A$ and some effective divisor $D$) and thus inside
some common zero-set of pluricanonical sections. We can actually
identify elements of ${\mathcal J}$ by looking at pluricanonical
sections whose vanishing order is above the stable vanishing order
by some appropriate amount.

\bigbreak\noindent(1.9) {\it Too Strongly Formulated Condition.} The
inequality $(\dagger)$ is weaker than
$$
\left|s_{\mathcal
J}\right|^2\,\sum_{j=1}^{q_k}\left|s^{(k)}_j\right|^2\leq
C_{m,k,{\mathcal
J}}\sum_{j=1}^{q_m}\left|s^{(m)}_j\right|^{\frac{2k}{m}}\leqno{(\dagger)^\sharp}
$$
for all $m,k\geq m_{\mathcal J}$ which states that the common
vanishing order of all $m$-canonical sections raised to the power
$\frac{k}{m}$ is no more than the common vanishing order of all
$k$-canonical sections multiplied by local generators of ${\mathcal
J}$.  By Skoda's result on ideal generation this stronger inequality
$(\dagger)^\sharp$ would imply immediately the finite generation of
the canonical ring.  See the Appendix (A.5).  We do not use this
formulation, because it is too strongly formulated and is much more
difficult to verify.

\bigbreak\noindent{\bf\S 2.} {\sc Constructing and Decreasing
Discrepancy Subspace.}  The proof of the finite generation of the
canonical ring is done by first constructing the initial
codimension-one discrepancy subspace and then decreasing the
discrepancy subspace until the stable vanishing order is achieved
everywhere.  The decreasing of the discrepancy subspace is done by
imitating the construction of the initial codimension-one
discrepancy subspace and using the family of subvarieties associated
to a closed positive $(1,1)$-current which is motivated by the
intuition of getting strict positive lower bound for the current
along the normal directions of the subvarieties with the modified
restriction of the current to the subvariety being of the special
form (that is, in the second case of the dichotomy (A.10)).

\bigbreak\noindent(2.1) {\it Construction of Initial Codimension-One
Discrepancy Subspace.}  As the first step we now construct the
initial codimension-one discrepancy subspace. We do this by using
the technique of the global generation of the multiplier ideal sheaf
(Appendix (A.1)) and the decomposition of $K_X$ as a sum of an ample
${\mathbb Q}$-line bundle and an effective ${\mathbb Q}$-divisor
from the general type property of $X$.

\medbreak Let $A$ be an ample line bundle on $X$ which is ample
enough for the global generation of multiplier ideal sheaves (see
Appendix (A.1)).  We write $aK_X=A+D$, where $D$ is an effective
divisor in $X$ and $a$ is a positive integer. We use the metric
$$\frac{1}{\Phi^m\left|s_D\right|^2}$$ for the line bundle
$$mK_X+D=\left(m+a\right)K_X-A,$$ where $s_D$ is the canonical
section of the line bundle $D$ so that the divisor of $s_D$ is
precisely $D$.  Let ${\mathcal I}^{(m)}$ be the multiplier ideal
sheaf of the metric
$$\frac{1}{\Phi^m\left|s_D\right|^2}.$$   Then the multiplier ideal sheaf
${\mathcal I}^{(m)}$ is generated by elements of
$$\displaylines{\qquad\qquad\qquad\Gamma\left(X,{\mathcal I}^{(m)}\left(mK_X+D+A\right)\right)\hfill\cr\hfill=
\Gamma\left(X,{\mathcal
I}^{(m)}\left(\left(m+a\right)K_X\right)\right)\subset
\Gamma\left(X,\left(m+a\right)K_X\right).\qquad\cr}$$  From the
Lemma on sup norm domination of metric by generators of multiplier
ideal (see Appendix (A.6)) we conclude that
$$
\left|s^{(k)}_j\right|^{\frac{2m}{k}}\left|s_D\right|^2\leq
C_{k,j}\sum_{j=1}^{q_{m+a}}\left|s^{(m+a)}_j\right|^2 \leqno{(\%)}
$$
for $k\in{\mathbb N}$.  This shows that we can choose ${\mathcal J}$
to be the ideal sheaf generated by $s_D$ and choose $m_{\mathcal J}$
as $a+1$.

\bigbreak\noindent(2.2) {\it Decreasing Discrepancy Subspace by
Imitating the Argument of Constructing Initial Codimension-One
Discrepancy Subspace.}  We are going to decrease the discrepancy
subspace by imitating the argument of constructing the initial
codimension-one discrepancy subspace given in (2.1).  Let us
re-interpret and recast the argument of (2.1) so that we can more
easily explain how we imitate it and adapt it to decrease the
discrepancy subspace.

\medbreak Let $q\in{\mathbb N}$.  Suppose $A$ is a holomorphic line
bundle on $X$ which is sufficiently ample so that for any point $P$
of $X$ there exist
$$
\tau_1,\cdots,\tau_r\in\Gamma\left(X,A\right)
$$
with the property that the multiplier ideal sheaf of the metric
$$
\frac{1}{\sum_{j=1}^r\left|\tau_j\right|^2}
$$
is in a neighborhood of $P$ in $X$ equal to ${\mathcal B}$ with
${\mathcal B}\subset\left({\mathfrak m}_P\right)^q$ and $P$ being an
isolated zero of ${\mathcal B}$. In other words, we are able to get
elements of $\Gamma\left(X,A\right)$ which could give metrics whose
multiplier ideal sheaves have the desired properties of isolated
zeroes at prescribed points.  We now have $aK_X=A+D$.  We can
interpret it as $aK_X-D=A$, which means that we can get elements of
$\Gamma\left(X,{\mathcal J}\left(aK_X\right)\right)=\Gamma\left(X,
aK_X-D\right)$ which could give metrics whose multiplier ideal
sheaves have the desired properties of {\it additional} isolated
zeroes at prescribed points, where ${\mathcal J}$ is the multiplier
ideal sheaf for the metric $\frac{1}{\left|s_D\right|^2}$.  More
precisely, for a point $P$ in $X$ we can find elements
$\sigma_1,\cdots,\sigma_\ell$ of $\Gamma\left(X, {\mathcal
J}\left(aK_X\right)\right)$ such that the multiplier ideal sheaf of
the metric
$$
\frac{1}{\sum_{j=1}^\ell\left|\sigma_j\right|^2}
$$
in a neighborhood of $P$ is of the form ${\mathcal A}{\mathcal J}$
where ${\mathcal A}\subset\left({\mathfrak m}_P\right)^q$ and $P$ is
an isolated zero of ${\mathcal A}$. The ideal sheaf ${\mathcal A}$
is what we mean by {\it additional} isolated zeroes.  Note that
${\mathcal J}$ is generated by the element $s_Ds_A\in\Gamma\left(X,
aK_X\right)$.

\medbreak What we need is actually not the global sections
$\sigma_1,\cdots,\sigma_\ell$ but the metrics
$$
\frac{1}{\sum_{j=1}^\ell\left|\sigma_j\right|^2}
$$
for $aK_X$ whose multiplier ideal sheaves have the desired
properties of {\it additional} isolated zeroes at prescribed points.
Of course, the use of such a metric is to enable us to consider the
metric
$$
\frac{e^{-m\varphi}}{\sum_{j=1}^\ell\left|\sigma_j\right|^2}
$$
of $(m+a)K_X$ so that the desired properties of {\it additional}
isolated zeroes at prescribed points would enable us to conclude
that $D$ can be used as a discrepancy subspace.

\medbreak Another important observation is that we can use the same
argument if we have such a metric not for $aK_X$ but for a ${\mathbb
Q}$-line bundle $aK_X+\delta B$ for some fixed ample line bundle $B$
and for a sufficiently small positive rational number $\delta$.  The
reason is the following.  Let us denote by $h$ this metric for
$aK_X+\delta B$.  In the argument what matters is the multiplier
ideal sheaves and not the metrics themselves.  We can write
$K_X=\alpha B+E$ for some rational positive number $\alpha$ and some
effective ${\mathbb Q}$-divisor $E$ and we can replace the use of
the metric $e^{-m\varphi}$ of $mK_X$ by the use of the metric
$$
\frac{
e^{-\left(m-\frac{\delta}{\alpha}\right)\varphi}}{\left|s_E\right|^{\frac{2\delta}{\alpha}}}\leqno{(\&)}
$$
of $mK_X-\delta B$ to form the metric
$$
\frac{
h\,e^{-\left(m-\frac{\delta}{\alpha}\right)\varphi}}{\left|s_E\right|^{\frac{2\delta}{\alpha}}}
$$
for $\left(m+a\right)K_X$ by putting together the metric $(\&)$ of
$mK_X-\delta B$ and the metric $h$ of $aK+\delta B$.  This is the
argument of absorbing a sufficiently small ample summand by the
bigness of the canonical line bundle (which is also described in
Appendix (A.12) for convenience of referral in some other steps of
this overview).  This is the part which needs the general type
assumption of $X$.  It is the same argument of absorption of small
ample summand by the bigness of the canonical line bundle which
makes the proof of the plurigenera for the case of general type so
much easier, because of the avoidance of a laborious estimation
process in analysis. This finishes our re-interpreting and recasting
of the construction of the initial codimension-one discrepancy
subspace given in (2.1).

\medbreak We now continue with the process of decreasing the
discrepancy subspace.  For the step of decreasing the discrepancy
subspace, as observed in (1.5), after blowing up we can assume
without loss of generality that the discrepancy ideal sheaf
${\mathcal J}$ is the ideal sheaf of an effective ${\mathbb
Z}$-divisor $Y=\sum_j b_j Y_j$ whose components are nonsingular
hypersurfaces $\left\{Y_j\right\}_j$ in normal crossing. To make the
argument simpler to understand, we look at the special case of $Y$
being a single nonsingular hypersurface $Y$ of $X$ with multiplicity
$1$. Let $\gamma$ be the generic stable vanishing order across $Y$.
Because of the skipping of diophantine arguments, we assume $\gamma$
to be rational.  We will replace $K_X$ on $X$ by $\left(K_X-\gamma
Y\right)|_Y$ on $Y$.  From the way the discrepancy subspace is
defined, we can assume without loss of generality that there is some
$m_1$-canonical section $s^*$ of $X$ such that at every point of
some Zariski open subset of $Y$ the multiplicity of $s^*$ is
strictly greater than $m_1\gamma$.   For the general case of
$Y=\sum_j b_j Y_j$ a modification of the argument works to reduce
one of the positive coefficients $b_j$ by $1$. Now for this special
case we have to look at $Y$ instead of $X$. There are a number of
modifications needed from construction of the initial
codimension-one discrepancy subspace given in (2.1).

\medbreak First of all, the ${\mathbb Q}$-line bundle
$\left(K_X-\gamma Y\right)|_Y$ may not be big on $Y$.  Let
$L=\left(K_X-\gamma Y\right)|_Y$.  Let $\Theta_Y$ be the curvature
current on $Y$ constructed from the metric
$\frac{e^{-\varphi}}{\left|s_Y\right|^{2\gamma}}$ of $K_X-\gamma Y$
on $X$. Some modification is needed to define $\Theta_Y$ on $Y$ (see
Appendix (A.11)).  It may happen that we have the extreme case of
the closed positive $(1,1)$-current $\Theta_L$ being of the special
form (that is, in the second case of the dichotomy (A.10)). In
general, there exist some nonnegative integer $\kappa\leq n-1$ and a
complex submanifold $V$ of complex dimension $\kappa$ in $Y$ such
that
\begin{itemize}
\item[(i)] the restriction of $\Theta_L$ to $V$ dominates some
strictly positive smooth $(1,1)$-form $\omega_V$ on $V$ (which we
can assume without loss of generality to be closed also), and
\item[(ii)] there exists a holomorphic family
$\left\{W_s\right\}_{s\in S}$ of subvarieties $W_s$ of pure
dimension $n-1-\kappa$ such that the ``modified restrictions'' of
$\Theta_Y$ to each $W_s$ are of the special form.
\end{itemize}
To explain ``modified restriction'', for the sake of simplicity let
us consider the case $\kappa=1$ so that each $W_s$ is a hypersurface
in $Y$.  In this case the modified restriction of $\Theta_Y$ to
$W_s$ is the restriction of $\Theta_Y-\gamma_{W_s}W_s$ to $W_s$,
where $\gamma_{W_s}$ is the generic Lelong number of $\Theta_Y$ at
points of $W_s$ [Siu 1974].   For the case of a general $\kappa$,
when $W_s$ is cut out by branches of divisors of pluricanonical
sections, $W_s$ is a hypersurface in a subvariety one dimension
higher, cut out by the branches of a subset of such divisors, which
is in another subvariety one dimension higher, cut out by the
branches of another subset of such divisors, {\it et cetera} until
we get to $Y$ and the modified restriction is defined inductively.
To define ``modified restriction'', instead of the divisor of a
pluricanonical section, we also allow ourselves the use of a
subspace defined by a metric of a pluricanonical line bundle with
nonnegative curvature current. For the general case the assumption
about ``modified restrictions'' of $\Theta_Y$ to $W_s$ means that
any ``modified restriction'' of $\Theta_Y$ to $W_s$ so defined must
be of the special form (that is, in the second case of the
dichotomy).  Since our purpose here is just to explain the main
argument, in order to avoid non-illuminating very complicated
descriptions and notations, so far as the assumption of the
``modified restriction'' of $\Theta_Y$ to $W_s$ is concerned, we
will confine ourselves to the case of $\kappa=1$.

\medbreak For $Y$ the subvariety $W_s$ will play the role of a point
$P$ in $X$ in our argument to decrease the discrepancy subspace. The
complex submanifold $V$ of $Y$ is related to the parameter space $S$
for the holomorphic family $\left\{W_s\right\}_{s\in S}$ of
subvarieties in $Y$ in that there is a holomorphic finite-fibered
map from $V$ to $S$. Intuitively, $V$ gives the directions in $Y$
where the current $\Theta_Y$ has a strictly positive lower bound and
each $W_s$ gives the directions where there is no longer any
positive lower bound for $\Theta_Y$.  The complex submanifold $V$
and the family $\left\{W_s\right\}_{s\in S}$ are constructed from
taking a fixed sufficiently ample line bundle $B$ on $X$ and using
elements of $\Gamma\left(Y, mL+B\right)$ for $m$ sufficiently large
so that
\begin{itemize}\item[(i)] the dimension of the restriction
$\Gamma\left(Y, mL+B\right)|V$ of $\Gamma\left(Y, mL+B\right)$ to
$V$ is of the order $cm^\kappa$ for some positive constant $c>0$ as
$m\to\infty$,
\item[(ii)] $W_s$ is obtained as the limit as $m\to\infty$ of the
multiplier ideal sheaves defined by appropriate roots of sums of
squares of elements of $\Gamma\left(Y, mL+B\right)|V$ which vanish
at a prescribed point to a high order depending on $m$.\end{itemize}
For more details about $W$ and $\left\{W_s\right\}_{s\in S}$ see
Appendix (A.13).

\medbreak Instead of using $D$ in $aK_X=A+D$ for $X$, we construct
our new divisor $D_Y$ in $Y$ by using the zero-set $\tilde D_Y$ of
the multiplier ideal sheaf of the closed positive $(1,1)$-current
$\left(\Theta_Y|_V\right)-\omega_V$ on $V$ (see Appendix (A.7)) and
the set of all $W_s$ such that $W_s$ intersects $\tilde D_Y$.  We
end up with a coherent ideal sheaf ${\mathcal K}$ on $X$ which
contains the ideal sheaf ${\mathcal I}_Y$ of $Y$ such that for each
$s\in S$ the metric of $\frac{m}{p_m}L+\frac{1}{p_m}B$ of the form
$$
\frac{1}{\sum_{j=1}^\ell\left|\sigma_j\right|^{\frac{2}{p_m}}}
$$
for some appropriate
$$
\sigma_1,\cdots,\sigma_\ell\in\Gamma\left(Y,{\mathcal
K}\left(mL+B\right)\right)
$$
(and some $m\in{\mathbb N}$ and $p_m\in{\mathbb N}$ with $p_m$
sufficiently large for the purpose of applying the technique in
Appendix (A.9)) gives the multiplier ideal sheaf ${\mathcal
K}{\mathcal I}_{W_s}$ with {\it additional} zeroes given by
${\mathcal I}_{W_s}$.  Now regard $\sigma_1,\cdots,\sigma_\ell$ as
elements of $\Gamma\left(Y,{\mathcal
I}_{m\varphi_Y}\left(mL+B\right)\right)$ (where
$e^{-\varphi_Y}=\frac{e^{-\varphi}}{\left|s_Y\right|^{2\gamma}}$)
and extend them to
$$\hat\sigma_1,\cdots,\hat\sigma_\ell\in\Gamma\left(X,{\mathcal
I}_{m\varphi_Y}\left(mL+B\right)\right).$$ This is possible by
Appendix (A.9) when $B$ is sufficiently ample (which we assume to be
the case).

\medbreak Also the ideal sheaf ${\mathcal K}$ can be obtained as the
multiplier ideal sheaf of a metric $\frac{\hat m}{\hat
p_m}L+\frac{1}{\hat p_m}B$ of the form
$$
\frac{1}{\sum_{j=1}^{\hat\ell}\left|\tau_j\right|^{\frac{2}{\hat
p_m}}}
$$
for some appropriate
$$
\tau_1,\cdots,\tau_{\hat\ell}\in\Gamma\left(Y,{\mathcal I}_{\hat
m\varphi_Y}\left(\hat mL+B\right)\right)
$$
(and some $\hat m\in{\mathbb N}$ and $\hat p_m\in{\mathbb N}$ with
$\hat p_m$ sufficiently large for the purpose of apply the technique
of in Appendix (A.9)).  Now extend $\tau_1,\cdots,\tau_{\hat\ell}$
to
$$\hat\tau_1,\cdots,\hat\tau_{\hat\ell}\in\Gamma\left(X,{\mathcal
I}_{\hat m\varphi_Y}\left(mL+B\right)\right).$$  Again this is
possible by Appendix (A.9) when $B$ is sufficiently ample (which we
assume to be the case).  Note that
$\hat\tau_1,\cdots,\hat\tau_{\hat\ell}$ can be regarded as the
analog of $s_Ds_A\in\Gamma\left(X, aK_X\right)$ in the construction
of the initial codimension-one discrepancy subspace done in (2.1)
and re-interpreted and recast above.  Note that whenever the
zero-set of ${\mathcal K}$ intersects some $W_s$, it contains all of
$W_s$, because of the way ${\mathcal K}$ is constructed from $\tilde
D_Y$ and the set of all $W_s$ such that $W_s$ intersects $\tilde
D_Y$.

\medbreak We now use the technique of constructing sections from
curvature currents of the special form on the subvariety $W_s$ and
extending them to all of $X$ by vanishing theorems from multiplier
ideal sheaves of the metrics described above.  We conclude that
${\mathcal K}{\mathcal I}_{m\varphi_Y}$ is generated by elements of
$\Gamma\left(X, mL\right)$ for all $m$ sufficiently large. Now we
use Appendix (A.6) to conclude that ${\mathcal K}$ is a discrepancy
ideal sheaf which is strictly bigger than ${\mathcal J}$.

\medbreak There are a number of subtle points which we gloss over in
our overview of the argument to decrease the discrepancy subspace.
Since it is easier to explain those subtle points in isolation, we
will do some of them in Appendix (A.14) and (A.15).

\medbreak\noindent(2.2.1) {\it Remark.} The above argument of
decreasing the discrepancy subspace by using the special form of the
current current and extension is a more streamlined version of the
argument in [Siu 2006] of slicing by ample divisors into curves and
using the second case of the dichotomy.  Instead of using ample
divisors to slice, we simply use a metric of the ample divisor with
nonnegative curvature current whose multiplier ideal sheaf gives the
ideal sheaf of a subvariety where the modified restriction of the
curvature current is of special form.  Special form here means the
second case of the dichotomy in [Siu 2006].

\bigbreak\noindent(2.3) {\it Abundance Conjecture.}  For a compact
complex algebraic manifold $W$ and an ample line bundle $A$ on $W$,
the abundance conjecture which compares
$$
\limsup_{m\to\infty}\frac{\log\dim_{\mathbb C}\Gamma\left(W,
mK_W+A\right))}{\log m}
$$
and
$$
\limsup_{m\to\infty}\frac{\log\dim_{\mathbb C}\Gamma\left(W,
mK_W\right)}{\log m}.
$$
deals with a situation similar to what is done in the argument (2.2)
to decrease the discrepancy subspace by using a holomorphic family
of subvarieties identified by directions of strict positive lower
bound for a closed positive $(1,1)$-current.  Unfortunately the
argument of absorption of small ample line bundle discussed in
(A.12) is used here, making it impossible to use directly the
argument for a proof of the abundance conjecture.  To adapt the
argument for use in a proof of the abundance conjecture, we
encounter the situation similar to adapting the proof of the
deformational invariance of the plurigenera for the case of general
type [Siu 1998] to the general algebraic case without the general
type assumption [Siu 2002], requiring an analytic argument of
controlling the estimates in passing to limit.

\bigbreak

\bigbreak\centerline{\bf APPENDIX}

\bigbreak\noindent(A.1) {\it Statement on Global Generation of
Multiplier Ideal Sheaves}. Let $L$ be a holomorphic line bundle over
an $n$-dimensional compact complex manifold $Y$ with a Hermitian
metric which is locally of the form $e^{-\xi}$ with $\xi$
plurisubharmonic. Let ${\cal I}_\xi$ be the multiplier ideal sheaf
of the Hermitian metric $e^{-\xi}$. Let $A$ be an ample holomorphic
line bundle over $Y$ such that for every point $P$ of $Y$ there are
a finite number of elements of $\Gamma(Y,A)$ which all vanish to
order at least $n+1$ at $P$ and which do not simultaneously vanish
outside $P$. Then $\Gamma(Y,{\cal I}_\xi\otimes(L+A+K_Y))$ generates
${\cal I}_\xi\otimes(L+A+K_Y)$ at every point of $Y$.

\bigbreak\noindent(A.2) {\it Skoda's Result on Ideal Generation}
[Skoda 1972]. Let $\Omega $ be a domain spread over $\mathbb{C}^{n}$
which is Stein. Let $\psi$ be a plurisubharmonic function on $\Omega
$, $g_{1},\ldots ,g_{p}$ be holomorphic functions on $\Omega $,
$\alpha
>1$, $q=\min \left(n,p-1\right)$, and $f$ be a holomorphic function on $\Omega
$. Assume that
$$\int_{\Omega }\frac{\left\vert f\right|^{2}e^{-\psi }}{\left(
\sum_{j=1}^p\left| g_{j}\right| ^{2}\right)^{\alpha q+1}}<\infty.$$
Then there exist holomorphic functions $h_{1},\ldots ,h_{p}$ on
$\Omega$ with $f=\sum_{j=1}^{p} h_jg_j$ on $\Omega$ such that
$$\int_{\Omega }\frac{\left\vert h_{k}\right|^{2}e^{-\psi
}}{\left(\sum_{j=1}^p\left|g_{j}\right|^2\right)^{\alpha q}}\leq
\frac{\alpha }{\alpha -1}\int_{\Omega
}\frac{\left|f\right|^2e^{-\psi }}{\left( \sum_{j=1}^p\left|
g_{j}\right|^2\right)^{\alpha q+1}}$$ for $1\leq k\leq p$.

\bigbreak\noindent(A.3) {\it Multiplier-Ideal Version of Skoda's
Result on Ideal Generation.} Let $X$ be a compact complex algebraic
manifold of complex dimension $n$, $L$ be a holomorphic line bundle
over $X$, and $E$ be a holomorphic line bundle on $X$ with metric
$e^{-\psi }$ such that $\psi $ is plurisubharmonic. Let $k\geq 1$ be
an integer, $G_{1},\ldots ,G_{p}\in \Gamma (X,L)$, and
$\left|G\right|^{2}=\underset{j=1}{\overset{p}{\sum
}}\left|G_j\right|^{2}$. Let
$\mathcal{I=I}_{(n+k+1)\log\left|G\right|^{2}+\psi}$ and
$\mathcal{J=I}_{(n+k)\log\left|G\right|^2+\psi}$. Then
$$\displaylines{\qquad\qquad\Gamma\left(X,\mathcal{I}\otimes
\left((n+k+1)L+E+K_{X}\right)\right)\hfill\cr\hfill
=\underset{j=1}{\overset{p}{%
\sum }}G_{j}\,\Gamma \left(X,\mathcal{J}\otimes\left(
(n+k)L+E+K_{X}\right)\right).\qquad\qquad\cr}$$

\bigbreak\noindent(A.4) {\it Finite Generation of Canonical Ring
From Achievement of Stable Vanishing Order.} Suppose the stable
vanishing orders are achieved at every point of $X$ for some
$m_0\in{\mathbb N}$. Denote $\left(m_0\right)!$ by $m_1$. Then the
canonical ring
$$
\bigoplus_{m=1}^\infty\Gamma\left(X, mK_X\right)
$$
is generated by
$$
\bigoplus_{m=1}^{\left(n+2\right)m_1}\Gamma\left(X, mK_X\right)
$$
and hence is finitely generated by the finite set of elements
$$
\left\{s^{(m)}_j\right\}_{1\leq m\leq m_1,\,1\leq j\leq q_m}.
$$

\medbreak\noindent{\it Proof.}  Let $e^{-\varphi}=\frac{1}{\Phi}$.
For $m>(n+2)m_1$ and any $s\in\Gamma\left(X, mK_X\right)$ we have
$$
\int_X\frac{\left|s\right|^2e^{-\left(m-\left(n+2\right)m_1-1\right)\varphi}}
{\left(\sum_{j=1}^{q_{m_1}}\left|s^{\left(m_1\right)}_j\right|^2\right)^{n+2}}<\infty,
$$
because $\left|s\right|^2\leq\tilde C\Phi^m$ on $X$ for some $\tilde
C$. By Skoda's theorem on ideal generation ((A.2) and (A.3)) there
exist
$$
h_1,\cdots,h_{q_{m_1}}\in\Gamma\left(X,\left(m-m_1\right)K_X\right)
$$
such that $s=\sum_{j=1}^{q_{m_1}}h_j s^{\left(m_1\right)}_j$.   If
$m-\left(n+2\right)m_1$ is still greater than $\left(n+2\right)m_1$,
we can apply the argument to each $h_j$ instead of $s$ until we get
$$
h^{(j_1,\cdots,j_\nu)}_1,\cdots,h^{(j_1,\cdots,j_\nu)}_{q_{m_1}}\in\Gamma\left(X,\left(m-m_1\left(\nu+1\right)\right)K_X\right)
$$
for $1\leq j_1,\cdots,j_\nu\leq q_{m_1}$ with $0\leq\nu< N$, where
$N=\left\lfloor\frac{m}{m_1}\right\rfloor$, such that
$$
s=\sum_{1\leq j_1,\cdots, j_N\leq
q_{m_1}}h^{(j_1,\cdots,j_{N-1})}_{j_N}\prod_{\lambda=1}^Ns^{\left(m_1\right)}_{j_\lambda}.
$$
Q.E.D.

\bigbreak\noindent(A.5) {\it Too Strongly Formulated Version of
Discrepancy Subspace.}  One can introduce a more strongly formulated
version of the discrepancy subspace by using the inequality
$$
\left|s_{\mathcal
J}\right|^2\,\sum_{j=1}^{q_k}\left|s^{(k)}_j\right|^2\leq
C_{m,k,{\mathcal
J}}\sum_{j=1}^{q_m}\left|s^{(m)}_j\right|^{\frac{2k}{m}}\leqno{(\dagger)_{k,m}^\sharp}
$$
for all $m,k\geq m_{\mathcal J}$, which states that the common
vanishing order of all $m$-canonical sections raised to the power
$\frac{k}{m}$ is no more than the common vanishing order of all
$k$-canonical sections multiplied by local generators of ${\mathcal
J}$.  This notion turns out to be too strong for our purpose.  It
actually gives immediately the finite generation of the canonical
ring.  We now verify that the inequality $(\dagger)^\sharp_{k,m}$
for all $m,k\geq m_{\mathcal J}$ implies right away the finite
generation of the canonical ring by Skoda's result on ideal
generation.  Take any $\lambda\in{\mathbb N}$. From
$$
\left(\sum_{j=1}^{q_k}\left|s^{(k)}_j\right|^2\right)^{\lambda}\leq
C_{\lambda,k}^\&\,\sum_{j=1}^{q_{\lambda k}}\left|s^{(\lambda
k)}_j\right|^2
$$
and
$$
\left|s_{\mathcal J}\right|^2\leq C^*_k\left|s_{\mathcal
J}\right|^{\frac{2}{k}}
$$
and from taking the $k$-root of $(\dagger)^\sharp_{k\lambda,m}$ that
$$
\left|s_{\mathcal
J}\right|^2\left(\sum_{j=1}^{q_k}\left|s^{(k)}_j\right|^2\right)^{\frac{\lambda}{k}}\leq
C^*_k\left(C_{\lambda,k}^\& C_{m,\lambda k,{\mathcal
J}}\sum_{j=1}^{q_{m}} \left|s^{(m)}_j\right|^{\frac{2\lambda
k}{m}}\right)^{\frac{1}{k}},
$$
which implies that for fixed $m$ and $\lambda$ we can find a
positive constant $\check C_{m,\lambda}$ and some
$\left(m,\lambda\right)$-dependent rapidly decreasing sequence of
positive numbers $\varepsilon_\ell$ used in the definition of
$$
\Phi=\sum_{\ell=1}^\infty\varepsilon_\ell\sum_{j=1}^{q_\ell}\left|s^{(\ell)}_j\right|^{\frac{2}{\ell}}
$$
such that
$$
\left|s_{\mathcal J}\right|^2\Phi^\lambda\leq\check
C_{m,\lambda}\left(\sum_{j=1}^{q_{m}}
\left|s^{(m)}_j\right|^2\right)^\lambda.\leqno{(\&)}
$$
Fix $m_0\geq m_{\mathcal J}$. Choose $\ell\in{\mathbb N}$
sufficiently large so that $\left|s_{\mathcal
J}\right|^{\frac{-2}{\ell}}$ is locally integrable on $X$. Replacing
$\lambda$ by $\ell\lambda$ and taking the $\ell$-th root of $(\&)$,
we get
$$
\frac{\Phi^\lambda}{\left(\sum_{j=1}^{q_{m}}
\left|s^{(m)}_j\right|^2\right)^\lambda}\leq\frac{\left(\check
C_{m,\ell\lambda}\right)^{\frac{1}{\ell}}}{\left|s_{\mathcal
J}\right|^{\frac{2}{\ell}}}.\leqno{(\&)^\prime}
$$
This implies that, for $\nu\geq(n+2)m$,
$$
\displaylines{\int_X\frac{\left|s^{(\nu+1)}_k\right|^2}{\Phi^{\nu-(n+2)m}\left(\sum_{j=1}^{q_m}\left|s^{(m)}_j\right|^2\right)^{(n+2)m}}
\leq
\int_X\frac{\left|s^{(\nu+1)}_k\right|^2}{\Phi^{\nu+1}}\,\frac{\Phi^{(n+2)m+1}}
{\left(\sum_{j=1}^{q_m}\left|s^{(m)}_j\right|^2\right)^{(n+2)m}}\cr
\leq \frac{1}{\left(\varepsilon_{\nu+1}\right)^{\nu+1}}
\int_X\frac{\Phi^{(n+2)m+1}}{\left(\sum_{j=1}^{q_m}\left|s^{(m)}_j\right|^2\right)^{(n+2)m}}
\leq \int_X\frac{\left(\check C_{m,\ell
m(n+2)}\right)^{\frac{1}{\ell}}\Phi}{\left(\varepsilon_{\nu+1}\right)^{\nu+1}\left|s_{\mathcal
J}\right|^{\frac{2}{\ell}}}<\infty,\cr}
$$
because
$$
\left|s^{(\mu+1)}_k\right|^2\leq\frac{\Phi^{\mu+1}}{\left(\varepsilon_{\mu+1}\right)^{\mu+1}}.
$$
We now apply the Multiplier-Ideal Version of Skoda's Theorem on
Ideal Generation (A.3) with the following choices.
\begin{itemize}
\item[(i)] $k=1$, $L=mK_X$, $E=\left(\nu-(n+2)m\right)K_X$,

\item[(ii)] $G_j=s^{(m)}_j$ for $1\leq j\leq q_m$.

\item[(iii)] $e^{-\psi}=\frac{1}{\Phi^{\nu-(n+1)m}}$.

\end{itemize}
We get
$$
s^{(\nu+1)}_k=\sum_{j=1}^m\sigma_{k,j}^{(\nu+1-m)}s^{(m)}_j
$$
for some
$\sigma_{k,j}^{(\nu+1-m)}\in\Gamma\left(X,(\nu+1-m)K_X\right)$ for
$1\leq j\leq q_m$.  We can now apply the same argument to
$\sigma^{(\nu+1-m)}_{k,j}$ instead of $s^{(\nu+1)}_k$ and continue
inductively.  This would give the finite generation of the canonical
ring.

\bigbreak\noindent(A.6) {\it Lemma on Sup Norm Domination of Metric
by Generators of Multiplier Ideal.} Let $f_j$ be holomorphic
functions on some open neighborhood $U$ of the origin in ${\mathbb
C}^n$. Let $\varepsilon_j>0$ and $m_j\in{\mathbb N}$ so that
$$
\Psi=\sum_{j=1}^\infty\varepsilon_j\left|f_j\right|^{\frac{2}{m_j}}
$$
converges uniformly on compact subsets of $U$.  Let ${\mathcal J}$
be the multiplier ideal sheaf of the metric $\frac{1}{\Psi}$ and
$g_1,\cdots,g_\ell$ be holomorphic function germs on ${\mathbb C}^n$
at the origin such that the stalk of ${\mathcal J}$ at the origin is
generated by $g_1,\cdots,g_\ell$ over ${\mathcal O}_{{\mathbb
C}^n,0}$. Then there exists an open neighborhood $W$ of the origin
in ${\mathbb C}^n$ where $g_1,\cdots,g_\ell$ are defined and there
exists a positive constant $C_j$ such that
$$
\left|f_j\right|^{\frac{2}{m_j}}\leq
C_j\sum_{k=1}^\ell\left|g_k\right|^2
$$
on $W$.

\medbreak\noindent(A.6.1) {\it Remark.}  The geometric reason for
this lemma is that the minimum of the orders of the zeros of the
generators of a multiplier ideal ${\mathcal J}$ should be no more
than the order of the pole of the metric $\frac{1}{\Psi}$.  A proof,
for example, is given in [Demailly 1992].

\bigbreak\noindent(A.7) {\it Multiplier Ideal Sheaves of Closed
Positive (1,1)-Currents.}  For a plurisubharmonic function $\varphi$
on some open subset $\Omega$ of ${\mathbb C}^n$, the multiplier
ideal sheaf ${\mathcal I}_\varphi$ is defined as consisting of all
holomorphic function germs $f$ on $\Omega$ such that
$\left|f\right|^2e^{-\varphi}$ is locally integrable.

\medbreak When we have two metrics $e^{-\varphi}$ and $e^{-\psi}$,
in general we do not have the relation ${\mathcal
I}_\varphi{\mathcal I}_\psi={\mathcal I}_{\varphi+\psi}$.  In many
situations considered in this article the two ideal sheaves
${\mathcal I}_\varphi{\mathcal I}_\psi$ and ${\mathcal
I}_{\varphi+\psi}$ are very close.  To emphasize the closeness in
such situations, we also use the notation $\widehat{{\mathcal
I}_\varphi{\mathcal I}_\psi}$ to denote ${\mathcal
I}_{\varphi+\psi}$ or even go to the extreme of simply using the
inaccurate expression ${\mathcal I}_\varphi{\mathcal I}_\psi$.

\medbreak Let $\Theta$ be a closed positive $(1,1)$-current on
$\Omega$. Locally there is a potential function $\psi$ for $\Theta$
in the sense that
$$
\Theta=\frac{\sqrt{-1}}{2\pi}\partial\bar\partial\psi
$$
locally.  The multiplier ideal sheaf ${\mathcal I}_\Theta$ is
defined as consisting of all holomorphic function germs $f$ on
$\Omega$ such that $\left|f\right|^2e^{-\psi}$ is locally
integrable.  The ideal sheaf ${\mathcal I}_\Theta$ is independent of
the choice of the local potential function $\psi$, because if
$\tilde\psi$ is another local potential function for $\Theta$, then
the difference $\psi-\tilde\psi$ must be pluriharmonic and therefore
must be smooth.

\bigbreak\noindent(A.8) {\it Existence of Global Sections of Amply
Twisted Multiple and Metric of Nonnegative Curvature Current.} For a
line bundle $L$ over a compact complex manifold $X$, the following
two statements are equivalent.
\begin{itemize}\item[(a)] There exists a metric $e^{-\psi}$ along the fibers of
$L$ such that the curvature current of $e^{-\psi}$ is a closed
positive $(1,1)$-current.
\item[(b)] For any sufficiently ample line bundle $A$ (which depends
on $X$ but independent of $L$) there is a nonzero element of
$\Gamma\left(X, mL+A\right)$ for any $m\in{\mathbb N}$.
\end{itemize}
The condition of sufficient ampleness of $A$ in Condition (b) is
satisfied if $A-K_X$ is sufficiently ample for the global generation
of multiplier ideal sheaves on $X$ (Appendix (A.1)).  The
implication of Condition (b) by Condition (a) simply comes from the
fact that ${\mathcal
I}_{m\psi}\left(mL+K_X+\left(A-K_X\right)\right)$ is globally
generated over $X$ for any $m\in{\mathbb N}$ and, as a consequence,
the subspace of $\Gamma\left(X, {\mathcal
I}_{m\psi}\left(mL+K_X+\left(A-K_X\right)\right)\right)$ of
$\Gamma\left(X, mL+A\right)$ is nonzero.

\medbreak For the other direction we take any nonzero section
$\sigma_m\in\Gamma\left(X, mL+A\right)$ for $m\in{\mathbb N}$ and
consider the closed positive $(1,1)$-current
$$
\Theta_m=\frac{1}{m}\frac{\sqrt{-1}}{2\pi}\partial\bar\partial\log\left|\sigma_m\right|^2
$$
which represents the class of $L+\frac{1}{m}A$.  Let $\omega_A$ be
the positive curvature form of a smooth metric $h_A$ of the line
bundle $A$.  The total mass of $\Theta_m$ with respect to the
K\"ahler metric $\omega_A$ is given by
$$\int_X\Theta_m\wedge\omega_A^{n-1}=\left(L+\frac{1}{m}A\right)A^{n-1}
$$
which is uniformly bounded for all $m\in{\mathbb N}$.  We can select
a subsequence $m_\nu$ so that $\Theta_{m_\nu}$ is weakly convergent
to $\Theta_\infty$ as $\nu\to\infty$.  The closed positive
$(1,1)$-current $\Theta_\infty$ represents the class of $L$ and, as
a consequence, we can find a metric $e^{-\psi}$ of $L$ whose
curvature current is equal to $\Theta_\infty$.

\bigbreak\noindent(A.9) {\it Extension of Global Twisted Sections of
Multiplier Ideal Sheaves.}  Let $X$ be a compact complex algebraic
manifold and $Y$ be a nonsingular hypersurface of $X$.  Then there
exists an ample line bundle $A$ with the property that, for any
holomorphic line bundle $L$ on $X$ with metric $e^{-\psi}$ such that
$\psi$ is locally plurisubharmonic, the map
$$
\Xi: \Gamma\left(X, {\mathcal I}_\psi\left(L+A\right)\right)\to
\Gamma\left(Y, \left({\mathcal I}_\psi\left/{\mathcal
I}_\psi{\mathcal I}_Y\right.\right)\left(L+A\right)\right)
$$
is surjective.

\medbreak\noindent{\it Proof.}  Choose an ample holomorphic line
bundle $A_0$ on $X$ so that the ideal sheaf of $Y$ is generated by
elements $\sigma_1,\cdots,\sigma_k\in\Gamma\left(X,A_0\right)$.
Choose an ample line bundle $A_2$ on $X$ with smooth positively
curved metric $h_2$ such that $A_2-K_X$ is ample with smooth
positively curved metric $h_3$.  Let $A=A_1+A_2+A_3$.  Use the
metric
$$
\frac{e^{-\psi}h_2h_3}{\sum_{j=1}^k\left|\sigma_j\right|^2}
$$
of $L+A-K_X$ whose multiplier ideal sheaf is ${\mathcal
I}_\psi{\mathcal I}_Y$. Then the vanishing of $H^1\left(X,{\mathcal
I}_\psi{\mathcal I}_Y\left(L+A\right)\right)$ implies the
surjectivity of the map $\Xi$.  Q.E.D.

\medbreak\noindent(A.9.1) {\it Remark.} The important point is that
the line bundle $A$ on $X$ needs only to be sufficiently ample and
this sufficient ampleness depends only on $Y$ and is independent of
$L$ and $e^{-\psi}$. Note that the comparison of ${\mathcal
I}_\psi\left/{\mathcal I}_\psi{\mathcal I}_Y\right.$ and ${\mathcal
I}_{\psi|Y}$ can be made by using the extension theorem of
Ohsawa-Takesgoshi [Ohsawa-Takesgoshi 1987] if the local
plurisubharmonic function $\psi$ is not identically equal to
$-\infty$ on $Y$.

\bigbreak\noindent(A.10) {\it Canonical Decomposition of Closed
Positive (1,1)-Current.} Let $\Theta$ be a closed positive
$(1,1)$-current on a complex manifold $X$.  Then $\Theta$ admits a
unique decomposition of the following form
$$ \Theta=\sum_{j=1}^J\gamma_j\left[V_j\right]+R,
$$
where $\gamma_j>0$, $J\in{\mathbb N}\cup\left\{0,\infty\right\}$,
$V_j$ is a complex hypersurface and the Lelong number of the
remainder $R$ is zero outside a countable union of subvarieties of
codimension $\geq 2$ in $X$ [Siu 1974].  We consider the dichotomy
into two cases.  The first case is either $R\not=0$ or $J=\infty$.
The second case is both $R=0$ and $J$ is finite.  We also say that
the current $\Theta$ is of the special form if it is in the second
case of the dichotomy.

\bigbreak\noindent(A.11) {\it Modified Restriction of Closed
Positive (1,1)-Current.}  Let $X$ be a compact complex algebraic
manifold of general type and let $e^{-\varphi}=\frac{1}{\Phi}$ be
the metric of minimum singularity as defined in (1.1).  Let
$$
\Theta_\varphi=\frac{\sqrt{-1}}{2\pi}\partial\bar\partial\varphi
$$
be the curvature current of the metric $e^{-\varphi}$ of $K_X$.  Let
$Y$ be a nonsingular hypersurface in $X$.  The generic stable
vanishing order $\gamma$ across $Y$ means the Lelong number of
$\Theta_\varphi$ at a generic point of $Y$ [Siu 1974].  Because of
the skipping of diophantine arguments, we assume $\gamma$ to be
rational.  The number $\gamma$ is also the infimum of the generic
vanishing order of $\left(s^{(m)}_j\right)^{\frac{1}{m}}$ across $Y$
for $m\in{\mathbb N}$ and $1\leq j\leq q_m$.

\medbreak We consider the question of how to define the restriction
to $Y$ of the closed positive $(1,1)$-current
$\Theta-\gamma\left[Y\right]$.  We call such a restriction to $Y$
the {\it modified restriction} of $\Theta$ to $Y$, because we are
restricting $\Theta$ after we modify it by subtracting
$\gamma\left[Y\right]$ from it.  Let $s_Y$ be the canonical section
of the line bundle associated to $Y$ so that the divisor of $s_Y$ is
precisely $Y$.  We know that $s^{(m)}_j$ vanishes to order at least
$m\gamma$ across $Y$ for each $m\in{\mathbb N}$ and each $1\leq
j\leq q_m$, but the multi-valued fraction
$$
\frac{\left(s^{(m)}_j\right)^{\frac{1}{m}}}{s_Y^\gamma}
$$
may still be identically zero on $Y$ for each $m\in{\mathbb N}$ and
each $1\leq j\leq q_m$ so that the sum
$$
\sum_{m=1}^\infty\varepsilon_m\sum_{j=1}^{q_m}
\left|\frac{\left(s^{(m)}_j\right)^{\frac{1}{m}}}{s_Y^\gamma}\right|^2
$$
may be identically zero on $Y$, making it impossible to consider the
metric
$$
\frac{1}{\sum_{m=1}^\infty\varepsilon_m\sum_{j=1}^{q_m}
\left|\frac{\left(s^{(m)}_j\right)^{\frac{1}{m}}}{s_Y^\gamma}\right|^2}
$$
of the ${\mathbb Q}$-line bundle $\left(L-\gamma Y\right)|_Y$ and to
use its curvature current as the restriction of the closed positive
$(1,1)$-current $\Theta-\gamma\left[Y\right]$ to $Y$.  The following
modification is needed in the process of constructing the
restriction of the closed positive $(1,1)$-current
$\Theta-\gamma\left[Y\right]$ to $Y$.  For $k\in{\mathbb N}$ let
$\gamma_k$ be the infimum of the vanishing order of the multi-valued
section $\left(s^{(m)}_j\right)^{\frac{1}{m}}$ across $Y$ for $1\leq
m\leq k$ and $1\leq j\leq q_m$. Consider the metric
$$
\frac{1}{\sum_{m=1}^k\varepsilon_m\sum_{j=1}^{q_m}
\left|\frac{\left(s^{(m)}_j\right)^{\frac{1}{m}}}{s_Y^{\gamma_k}}\right|^2}
$$
of the ${\mathbb Q}$-line bundle $\left.\left(L-\gamma_k
Y\right)\right|_Y$ on $Y$ and its curvature current
$$
\Theta_k=\frac{\sqrt{-1}}{2\pi}\partial\bar\partial\log\sum_{m=1}^k\varepsilon_m\sum_{j=1}^{q_m}
\left|\frac{\left(s^{(m)}_j\right)^{\frac{1}{m}}}{s_Y^{\gamma_k}}\right|^2
$$ which is a
closed positive $(1,1)$-current on $Y$.  We know that the sequence
$\gamma_k$ is non-increasing and its limit is $\gamma$ as
$k\to\infty$ so that the ${\mathbb Q}$-line bundle
$\left.\left(L-\gamma_k Y\right)\right|_Y$ on $Y$ approaches the
${\mathbb Q}$-line bundle $\left.\left(L-\gamma Y\right)\right|_Y$
on $Y$ as $k\to\infty$.  The restriction of the closed positive
$(1,1)$-current $\Theta-\gamma\left[Y\right]$ to $Y$ can be defined
as the (weak) limit of $\Theta_k$ (or its subsequence) as
$k\to\infty$.

\bigbreak\noindent(A.12) {\it Absorption of Small Ample Line Bundle
and Small Modification of Metric.} Let $X$ be a compact complex
algebraic manifold of general type. Let $B$ be an ample line bundle
on $X$.  Then $K_X=\alpha B+E$ for some rational positive number
$\alpha$ and some effective ${\mathbb Q}$-divisor $E$.  Let
$e^{-\varphi}$ be the metric of minimum singularity for $K_X$.
Consider the metric
$$
h_\delta:=\frac{
e^{-\left(m-\frac{\delta}{\alpha}\right)\varphi}}{\left|s_E\right|^{\frac{2\delta}{\alpha}}}
$$
of $mK_X-\delta B$.  Let $\delta>0$ and $p\in{\mathbb N}$ and
$e^{-\psi}$ be a metric of nonnegative curvature current for
$pK_X+\delta B$.   We form the metric $h_\delta e^{-\psi}$ of
$(p+m)K_X$ and describe it as the {\it absorption of the small ample
line bundle} $\delta B$.   Let $h_B$ be a strictly positively curved
smooth metric of $B$.  We also call the metric
$h_\delta\left(h_B\right)^\delta$ of $mK_X$ a {\it small
modification} of the metric $e^{-m\varphi}$ of $mK_X$.

\bigbreak\noindent(A.13) {\it Family of Subvarieties Associated to a
Closed Positive (1,1)-Current.}

\bigbreak Let $X$ be a compact complex algebraic manifold and
$\Theta$ be a closed positive $(1,1)$-current on $X$.  We are going
to associate to $\Theta$ a family of subvarieties.  The motivation
is that there is a lower bound of positivity for $\Theta$ in the
normal directions of the subvarieties and the subvarieties are
minimum with respect to this property.  Let $\theta$ be the
$(1,1)$-class (an element of $H^1\left(X,\Omega^1_X\right)$) which
is defined by $\Theta$.  The class $\theta$ may not be an integral
class (that is, it may not come from $H^2\left(X, {\mathbb
Z}\right)$).

\medbreak By the simultaneous approximation of a finite collection
of real numbers by rational numbers (see, for example, [Hardy-Wright
1960, p.170, Th.200]) we can find elements $\phi_m$ of
$H^1\left(X,\Omega^1_X\right)$ and $\delta>0$ such that
$m\left(\theta+\phi_m\right)$ comes from $H^2\left(X, {\mathbb
Z}\right)$ and $m^{1+\delta}\phi_m\to 0$ in
$H^1\left(X,\Omega^1_X\right)$ as $m\to\infty$.  Let $L_m$ be the
holomorphic line bundle on $X$ which corresponds to the integral
$(1,1)$-class $m\left(\theta+\phi_m\right)$.  We choose an ample
line bundle $A_1$.  Since $\delta$ is positive and
$m^{1+\delta}\phi_m\to 0$ as $m\to\infty$, there exists
$m_0\in{\mathbb N}$ such that for $m\geq m_0$ we can find
$q_m\in{\mathbb N}$ and a metric $e^{-\varphi_m}$ for the ${\mathbb
Q}$-line bundle $L_m+\frac{q_m}{m}A_1$ such that $\varphi_m$ is
locally plurisubharmonic and its curvature current $\Theta_m$
approaches $\Theta$ and $\frac{q_m}{m}\to 0$ as $m\to\infty$.

\medbreak Let $A$ be a holomorphic line bundle on $X$ such that
$A-K_X$ is ample enough for the global generation of multiplier
ideal sheaves on $X$ (Appendix (A.1)).  Define the number $\kappa$
by
$$
\limsup_{m\to\infty}\frac{\log\dim_{\mathbb C}\Gamma\left(X,
{\mathcal I}_{\varphi_m}\left(L_m+A\right)\right)}{\log m}.
$$
This is a way of measuring the lower bound of positivity, because we
have global generation of the multiplier ideal sheaf when we add $A$
to $L_m$ to consider $\Gamma\left(X, {\mathcal
I}_{\varphi_m}\left(L_m+A\right)\right)$.

\medbreak Suppose in a neighborhood $U$ of a point $P$ of $X$ the
closed positive $(1,1)$-current $\Theta$ dominates a strictly
positive smooth $(1,1)$-current $\omega$ for $m$ sufficiently large.
Then $\Theta_m$ dominates $\frac{m}{2}\omega$.  We can then use a
local coordinate system $z=\left(z_1,\cdots,z_n\right)$ centered at
$P$ so that the unit ball $B$ is relatively compact in $U$.  We can
consider the function $\tilde\chi$ on $X$ which is equal to
$\frac{1}{\left|z\right|^2}$ on $B$ and equal to $1$ on $X-B$ and we
smooth $\chi$ slightly near the boundary of $B$ to get $\chi$ which
is smooth and positive on $X-\left\{P\right\}$ and is equal to
$\tilde\chi$ outside a very small neighborhood of the boundary of
$B$ in $U$.  Then for some $p_0\in{\mathbb N}$ sufficiently large
$$\frac{p_0}{2}\omega+\frac{\sqrt{-1}}{2\pi}\partial\bar\partial\log\chi$$
is a positive $(1,1)$-current on $U$.  Thus the curvature current of
the metric $e^{-\varphi_m}\chi^q$ is a positive $(1,1)$-current on
all of $X$ for $m\geq p_0q$.

\medbreak By Skoda's result (A.2), we know that the multiplier ideal
sheaf of $e^{-\varphi_m}\chi^q$ is contained in
$$
\left({\mathfrak m}_P\right)^{q-n-1}{\mathcal I}_{\varphi_m}
$$
for $m\geq p_0q$.  From the surjectivity of
$$
\Gamma\left(X,{\mathcal I}_{\varphi_m}\left(L_m+A\right)\right)\to
{\mathcal I}_{\varphi_m}\left/\left(\left({\mathfrak
m}_P\right)^{q-n-1}{\mathcal I}_{\varphi_m}\right)\right.
$$
for $m\geq p_0q$ we conclude that
$$
\dim_{\mathbb C}\Gamma\left(X,{\mathcal
I}_{\varphi_m}\left(L_m+A\right)\right)\geq\dim_{\mathbb
C}\left({\mathcal O}_X\left/\left({\mathfrak
m}_P\right)^{q-n-1}\right.\right)={q-2\choose n}
$$
for $m\geq p_0q$ and
$$
\dim_{\mathbb C}\Gamma\left(X,{\mathcal
I}_{\varphi_m}\left(L_m+A\right)\right)\geq{\left\lfloor
\frac{m}{p_0}\right\rfloor-2\choose n}.
$$
This implies that
$$
\limsup_{m\to\infty}\frac{\log\dim_{\mathbb C}\Gamma\left(X,
{\mathcal I}_{\varphi_m}\left(L_m+A\right)\right)}{\log m}\geq n
$$
and $\kappa\geq n$.

\medbreak We now consider another case.  Suppose $V$ is a complex
submanifold of complex dimension $d$ in $X$ and $\Theta|_V$
dominates some strictly positive smooth $(1,1)$-current $\omega_V$
on the neighborhood of some point $P$ in $V$.  Assume that $A$ is
sufficiently ample so that the map
$$
\Gamma\left(X,{\mathcal I}_{\varphi_m}\left(L_m+A\right)\right)\to
\Gamma\left(V,{\mathcal I}_{\varphi_m|_V}\left(L_m+A\right)\right)
$$
is surjective (A.9).  Since our earlier argument gives
$$
\dim_{\mathbb C}\Gamma\left(V,{\mathcal
I}_{\varphi_m|_V}\left(L_m+A\right)\right)\geq cm^d
$$
for some $c>0$ and for all $m\in{\mathbb N}$ sufficiently large, it
follows that
$$
\limsup_{m\to\infty}\frac{\log\dim_{\mathbb C}\Gamma\left(X,
{\mathcal I}_{\varphi_m}\left(L_m+A\right)\right)}{\log m}\geq d
$$
and $\kappa\geq d$.

\medbreak Conversely, suppose we have the growth rate and we would
like to conclude about the positive lower bound of the curvature
current when restricted to some submanifold.  We again choose some
transversal submanifold $V$ of complex dimension $\kappa$ which is
independent of $m$ in the sense that
$$
\dim_{\mathbb C}\Gamma\left(V,{\mathcal
I}_{\varphi_m|_V}\left(L_m+A\right)\right)\geq cm^\kappa
$$
for some $c>0$ and for all $m\in{\mathbb N}$ sufficiently large. The
choice of $V$ simply means that it is transversal to proper
subvarieties of $X$ defined by homogeneous polynomials of elements
of $\Gamma\left(X,{\mathcal
I}_{\varphi_m}\left(L_m+A\right)\right)$.  Those subvarieties are
from a countable number of holomorphic families of subvarieties.  We
consider elements in $\Gamma\left(X,{\mathcal
I}_{\varphi_m}\left(L_m+A\right)\right)$ whose restrictions to $V$
vanish to order $q$ at a point $P$ of $V$.  We now consider their
common zero-set $W$.  If the vanishing order at $P$ or any other
points of $W$ is substantially less than $\beta q$ for some positive
$\beta$ sufficiently close to $0$, for example $\beta\leq
q^{-\delta}$ for some $0<\delta<1$, then we can create a metric for
$L_m$ (by raising $m$ to a high multiple first and then taking root
of the metric later) so that the restriction to some $\tilde V$ with
complex dimension greater than that of $V$ gives a metric of
isolated high order singularity at a point of $\tilde V$.  This
would mean that the growth order is greater than $\kappa$ which is a
contradiction.  For every $P$ we can form the subvariety $W$ which
depends on $P$ and which we denote by $W_P$.  As $P$ varies inside
$V$, we can get a holomorphic family of subvarieties
$\left\{W_s\right\}_{s\in S}$.

\bigbreak\noindent(A.14) {\it Additional Vanishing and Minimal
Center of Log Canonical Singularity for the Second Case of
Dichotomy.}  Let $X$ be a compact complex algebraic manifold and $Y$
be a nonsingular complex hypersurface in $X$.  Let $\gamma$ be the
generic stable vanishing order across $Y$ which we assume to be
rational, because we are skipping the diophantine arguments.  Let
$\frac{1}{\Phi}=e^{-\varphi}$ be the metric of $K_X$ of minimum
singularity with curvature current $\Theta_\varphi$.  We have the
modified restriction $\Theta_\varphi-\gamma\left[Y\right]$ defined
in Appendix (A.11), which we denote by $\Theta_Y$.  We assume that
the current $\Theta_Y$ is of the special form ({\it i.e.,} in the
second case of the dichotomy according to the terminology of [Siu
2006]), which means that there are a finite number of complex
hypersurfaces $V_j$ in $Y$ for $1\leq j\leq J$ such that
$$
\Theta_Y=\sum_{j=1}^J\alpha_j\left[V_j\right].
$$
Again because we are skipping the diophantine argument, we assume
that each $\alpha_j$ is rational.  Take a sufficiently divisible
positive integer $\tilde m$ so that $\tilde m\gamma$ and $\tilde
m\alpha_j$ are positive integers.  We have a holomorphic section
$$
\sigma:=\prod_{j=1}^J\left(s_{V_j}\right)^{\tilde m\alpha_j}
$$
over $Y$ of the tensor product of the holomorphic line bundle
$\tilde m\left(K_X-\gamma Y\right)$ and some flat bundle over $Y$.
The flat bundle situation is handled by a modification of the
technique of Shokurov [Shokurov 1985] using the theorem of
Riemann-Roch and the vanishing theorem for multiplier ideal sheaves.
By Shokurov's technique we can get a section of $\tilde
m\left(K_X-\gamma Y\right)$ over $Y$ with the same divisor as
$\sigma$.  Here we are interested in explaining a technique of
minimum center of log canonical singularity and we will not go into
Shokurov's technique. We would like to extend $\sigma$ to an element
of $\Gamma\left(X,\tilde m\left(K_X-\gamma Y\right)\right)$.  The
extension is done by using the vanishing theorem for multiplier
ideal sheaves.

\medbreak Let $L=K_X-\gamma Y$.  If the generic stable vanishing
order $\gamma$ across $Y$ is achieved by some
$\tau\in\Gamma\left(X,\hat mK_X\right)$, then
$$
\left.\frac{\tau}{\left(s_Y\right)^{\hat m\gamma}}\right|_Y
$$
would be a holomorphic section of $\hat m L=\hat m\left(K_X-\gamma
Y\right)$ over $Y$ and its divisor must be $\hat
m\sum_{j=1}^J\alpha_j V_j$, because the curvature current $\Theta_Y$
is of the special form.  There is nothing more to do and there is no
need for any further discussion.

\medbreak What we are interested in is the case when the generic
stable vanishing order $\gamma$ across $Y$ is not achieved, which we
now assume to be the case.  For $N\in{\mathbb N}$ consider the
$N$-truncation
$$
\Phi_N=\sum_{m=1}^N\varepsilon_m\sum_{j=1}^{q_m}\left|s^{(m)}_j\right|^{\frac{2}{m}}.
$$
Then no matter how large $N$ is, the vanishing order
$\tilde\gamma_N$ of $\Phi_N$ across $Y$ is still strictly greater
than $\gamma$. Because the vanishing theorem adds one $K_X$, in
order to get $L-Y=K_Y-(1+\gamma)Y$ (which is to replace the new copy
of $K_X$ by $L$ and to get one order of vanishing across $Y$ needed
for the vanishing of the first cohomology to extend sections from
$Y$ to $X$) from $K_X$ we have to get an extra vanishing order
$1+\gamma$ across $Y$. For the extension of $\sigma$ to an element
of $\Gamma\left(X,\tilde m\left(K_X-\gamma Y\right)\right)$, the
ideal situation is to get the vanishing of $H^1\left(X,{\mathcal
I}_Y\left(\tilde m L\right)\right)$ so that we have the surjectivity
of the map
$$
\Gamma\left(X, \tilde m L\right)\to \Gamma\left(Y, \tilde m
L\right).
$$
The cohomology group $H^1\left(X,{\mathcal I}_Y\left(\tilde m
L\right)\right)$ is the same as $H^1\left(X,\left({\mathcal
I}_Y\right)^{\tilde m\gamma+1}\left(\tilde m K_X\right)\right)$.  In
order to get its vanishing by using multiplier ideal sheaf, because
of the addition of one $K_X$, the ideal situation is to get a metric
for $\left(\tilde m-1\right)K_X$ whose multiplier ideal sheaf is
$\left({\mathcal I}_Y\right)^{\tilde m\gamma+1}$.  If we use the
metric $\frac{1}{\Phi^{\tilde m-1}}=e^{-\left(\tilde
m-1\right)\tilde\varphi}$ for $\left(\tilde m-1\right)K_X$, the
multiplier ideal sheaf at a generic point of $Y$ is $\left({\mathcal
I}_Y\right)^{\left(\tilde m-1\right)\gamma}$.  Compared to the
required $\left({\mathcal I}_Y\right)^{\tilde m\gamma+1}$, we need
to add a factor of $\left({\mathcal I}_Y\right)^{\gamma+1}$ to
$\left({\mathcal I}_Y\right)^{\left(\tilde m-1\right)\gamma}$.  We
can achieve this at a generic point of $Y$ by using the interpolated
metric
$$
e^{-\psi_\delta}:=\left(\frac{1}{\Phi^{1-\delta}\Phi_N^\delta}\right)^{\tilde
m-1}
$$
for $\left(\tilde m-1\right)K_X$ with $0<\delta<1$ determined by
$$
\left(\tilde
m-1\right)\left((1-\delta)\gamma+\delta\tilde\gamma_N\right)=\tilde
m\gamma+1,
$$
which means
$$
\delta=\frac{1+\gamma}{\left(\tilde
m-1\right)\left(\tilde\gamma_N-\gamma\right)}.
$$
To get $0<\delta<1$, by replacing $\tilde m$ by a high multiple, we
can assume that $$\tilde m>\frac{\tilde\gamma_N-\gamma}{1+\gamma}.$$

\medbreak The main difficulty is that the multiplier ideal sheaf
${\mathcal I}_{\psi_\delta}$ of $e^{-\psi_\delta}$ may be smaller
than $\left({\mathcal I}_Y\right)^{\tilde m\gamma+1}$ so that we can
only get
$$
H^1\left(X, {\mathcal I}_{\psi_\delta}\left(\tilde m
L\right)\right)=0
$$
and the surjectivity of
$$
\Gamma\left(X, \tilde m L\right)\to \Gamma\left(X, \left({\mathcal
O}_X\left/{\mathcal I}_{\psi_\delta}\right.\right)\left(\tilde m
L\right)\right).
$$
Since $\sigma$ is only an element of $\Gamma\left(Y, \tilde m
L\right)$ and may not be an element of $\Gamma\left(X,
\left({\mathcal O}_X\left/{\mathcal
I}_{\psi_\delta}\right.\right)\left(\tilde m L\right)\right)$, in
general we would have trouble extending $\sigma$ to all of $X$.  We
call this difficulty {\it additional vanishing}, because the
vanishing order of
$\left(\Phi^{1-\delta}\Phi_N^\delta\right)^{\tilde m-1}$ is more
than the desired order of $\tilde m\gamma+1$ across $Y$ and there is
additional vanishing.  We do not even know that the zero-set of
${\mathcal I}_{\psi_\delta}$ is just $Y$.  There may even be zeroes
of ${\mathcal I}_{\psi_\delta}$ outside $Y$.

\medbreak Here we only want to explain one special technique which
handles the difficulty of additional vanishing for the case when the
zero-set of ${\mathcal I}_{\psi_\delta}$ is contained in $Y$. The
section $\sigma$ can be considered as an element of $$\Gamma\left(Y,
\left(\left({\mathcal I}_Y\right)^{\tilde
m\gamma}\left/\left({\mathcal I}_Y\right)^{\tilde
m\gamma+1}\right.\right)\left(\tilde mK_X\right)\right).$$  The
difficulty of additional vanishing means that ${\mathcal
I}_{\psi_\delta}$ does not contain $\left({\mathcal
I}_Y\right)^{\tilde m\gamma+1}$ and as a consequence $\sigma$ may
not be an element of $\Gamma\left(Y, \left({\mathcal
O}_X\left/{\mathcal I}_{\psi_\delta}\right.\right)\left(\tilde
mK_X\right)\right)$.  The technique is to decrease $\delta$ to
$\tilde\delta$ such that $\sigma$ can induce a well-defined and not
non-identically-zero element $\tilde\sigma$ of $\Gamma\left(Y,
\left({\mathcal O}_X\left/{\mathcal
I}_{\psi_{\tilde\delta}}\right.\right)\left(\tilde
mK_X\right)\right)$.  Then from $$H^1\left(X, {\mathcal
I}_{\psi_{\tilde\delta}}\left(\tilde m L\right)\right)=0
$$
and the surjectivity of
$$
\Gamma\left(X, \tilde m L\right)\to \Gamma\left(X, \left({\mathcal
O}_X\left/{\mathcal
I}_{\psi_{\tilde\delta}}\right.\right)\left(\tilde m L\right)\right)
$$
we would be able to extend $\tilde\sigma$ to an element of
$\hat\sigma$ of $\Gamma\left(X, \tilde m L\right)$.  Since the
current $\Theta_Y$ is of the special form and since $\hat\sigma$
restricts to $\tilde\sigma$ on $\left({\mathcal O}_X\left/{\mathcal
I}_{\psi_{\tilde\delta}}\right.\right)\left(\tilde m L\right)$, it
follows that $\hat\sigma$ must induce $\sigma$ as an element of
$$\Gamma\left(Y, \left(\left({\mathcal I}_Y\right)^{\tilde
m\gamma}\left/\left({\mathcal I}_Y\right)^{\tilde
m\gamma+1}\right.\right)\left(\tilde mK_X\right)\right).$$

\medbreak Let us explain the above procedure with the following
local concrete example.  Assume that $\gamma=0$ so that $\tilde
m\gamma=0$.  Locally at the origin $P_0$ of a local coordinate
system of $\left(z_1,\cdots,z_n\right)$ of $X$ we suppose that $Y$
is given by $z_1=0$ and $\Theta_Y$ is just the hypersurface $V_1$
given $z_2=0$ with coefficient $1$.  Assume that
$\psi_{\tilde\delta}$ is generated by $z_1^2$ and $z_1z_2$.  Near
$P_0$ the ringed space $\left(Y,{\mathcal O}_X\left/{\mathcal
I}_{\psi_{\tilde\delta}}\right.\right)$ is reduced except along $V$.
Assume that near $P_0$ the divisor $\sigma$ is locally the
hypersurface given by $z_2=0$ with coefficient $1$. When we consider
the element in ${\mathcal O}_X\left/{\mathcal
I}_{\psi_{\tilde\delta}}\right.$ induced by $\sigma$ locally near
$P_0$ we simply get a section of the conormal bundle of $V$ in $Y$
which is non-identically-zero (from differentiating $\sigma$ in the
normal direction of $V$).

\bigbreak\noindent(A.15) {\it Positive Lower Bound of Curvature
Current in Ambient Space.}  Let $X$ be a compact complex manifold of
general type with metric $\frac{1}{\Phi}$ of minimum singularity and
its $N$-th truncation $\frac{1}{\Phi_N}$.  Let $Y$ be a nonsingular
complex hypersurface in $X$ whose generic stable vanishing order
$\gamma$ is not achieved.  We choose two appropriate sufficiently
large positive integers $m$ and $N$ and construct two interpolations
$$
e^{-\psi_\nu}=\left(\frac{1}{\Phi^{1-\delta_\nu}\Phi_N^{\delta_\nu}}\right)^m\quad{\rm
for\ \ }\nu=1, 2
$$
with $0<\delta_1<\delta_2<1$ very close to each other and separated
by some critical value for integrability at a generic point of $Y$
(also after appropriate small modifications (A.12) if necessary) so
that the two multiplier ideal sheaves ${\mathcal I}_{\psi_1}$ and
${\mathcal I}_{\psi_2}$ satisfy the relation ${\mathcal
I}_{\psi_2}={\mathcal I}_Y{\mathcal I}_{\psi_1}$.  From $H^p\left(X,
{\mathcal I}_{\psi_\nu}\left((m+1)K_X\right)\right)=0$ for $p\geq 1$
it follows that
$$
H^p\left(Y, {\mathcal I}_{\psi_1}{\mathcal
O}_Y\left((m+1)K_X\right)\right)=0\quad{\rm for\ \ }p\geq 1,
$$
where ${\mathcal I}_{\psi_1}{\mathcal O}_Y={\mathcal
I}_{\psi_1}\left/\left({\mathcal I}_Y{\mathcal
I}_{\psi_1}\right)\right.$.  When $E$ is a holomorphic line bundle
on $X$ with metric $e^{-\chi}$ of nonnegative curvature current, we
can also conclude in the same way after a careful choice of
$\delta_1$ and $\delta_2$ and appropriate small modifications that
$$
H^p\left(Y, \widehat{{\mathcal I}_\chi{\mathcal
I}_{\psi_1}}{\mathcal O}_Y\left(E+(m+1)K_X\right)\right)=0\quad{\rm
for\ \ }p\geq 1.
$$
(See (A.7) for the meaning of the symbol $\widehat{{\mathcal
I}_\chi{\mathcal I}_{\psi_1}}$.) The point is that we get a
vanishing result on $Y$ from vanishing results on $X$.  The
vanishing theorem for multiplier ideal sheaves requires strictly
positive lower bound for the curvature current. Sometimes we have
such strictly positive lower bound for its curvature current on $X$,
but not on the restriction to $Y$.  In such a case this method
enables us to get a vanishing result by using the strictly positive
lower bound for the curvature current from the ambient space $X$.
For example, though $X$ is of general type, $Y$ may not be of
general type.  Also this is another way to handle additional
vanishing discussed in (A.14) so that we deal with the reduced
structure of $Y$ instead of an unreduced structure for $Y$ (which
comes from the additional vanishing).

\bigbreak\noindent{\it References}

\medbreak\noindent[Birkan-Cascini-Hacon-McKernan 2006] C. Birkar, P.
Cascini, C. Hacon, and J. McKernan, Existence of minimal models for
varieties of log general type, arXiv:math/0610203.

\medbreak\noindent[Demailly 1992] J.-P. Demailly, Regularization of
closed positive currents and Intersection Theory, {\it J. Alg.
Geom.} {\bf 1} (1992) 361-409.

\medbreak\noindent[Hardy-Wright 1960] G. H. Hardy and E. M. Wright,
{\it An Introduction to the Theory of Numbers}, 4th ed., Oxford
University Press 1960.

\medbreak\noindent[Kawamata 1982] Y. Kawamata, A generalization of
Kodaira-Ramanujam's vanishing theorem. {\it Math. Ann.} {\bf 261}
(1982), 43-46.

\medbreak\noindent[Kawamata 1985] Y. Kawamata, Pluricanonical
systems on minimal algebraic varieties. {\it Invent. Math.} {\bf 79}
(1985), 567--588.

\medbreak\noindent[Nadel 1990] A. Nadel, Multiplier ideal sheaves
and K\"ahler-Einstein metrics of positive scalar curvature. {\it
Ann. of Math.} {\bf 132} (1990), 549-596.

\medbreak\noindent[Ohsawa-Takegoshi 1987] Ohsawa, T., Takegoshi, K.:
On the extension of $L^2$ holomorphic functions, {\it Math.
Zeitschr.} {\bf 195}, 197-204 (1987).

\medbreak\noindent[Paun 2005] M. Paun, Siu's invariance of
plurigenera: a one-tower proof, {\it J. Differential Geom.} {\bf 76}
(2007), no. 3, 485--493.

\medbreak\noindent[Shokurov 1985] V.~V. Shokurov, A nonvanishing
theorem. {\it Izv. Akad. Nauk SSSR Ser. Mat.} {\bf 49} (1985),
635--651.

\medbreak\noindent[Siu 1974] Y.-T. Siu, Analyticity of sets
associated to Lelong numbers and the extension of closed positive
currents. {\it Invent. Math.} {\bf 27} (1974), 53-156.

\medbreak\noindent[Siu 1996] Siu, Y.-T., The Fujita conjecture and
the extension theorem of Ohsawa-Takegoshi, in {\it Geometric Complex
Analysis} ed. Junjiro Noguchi {\it et al}, World Scientific:
Singapore, New Jersey, London, Hong Kong 1996, pp. 577-592.

\medbreak\noindent[Siu 1998]  Y.-T. Siu, Invariance of plurigenera,
{\it Invent. Math.} {\bf 134} (1998), 661-673.

\medbreak\noindent[Siu 2002]  Y.-T. Siu, Extension of twisted
pluricanonical sections with plurisubharmonic weight and invariance
of semipositively twisted plurigenera for manifolds not necessarily
of general type. In: {\it Complex Geometry: Collection of Papers
Dedicated to Professor Hans Grauert}, Springer-Verlag 2002,
pp.223-277.

\medbreak\noindent[Siu 2003] Y.-T. Siu, Invariance of plurigenera
and torsion-freeness of direct image sheaves of pluricanonical
bundles, In: {\it Finite or Infinite Dimensional Complex Analysis
and Applications} (Proceedings of the 9th International Conference
on Finite or Infinite Dimensional Complex Analysis and Applications,
Hanoi, 2001), edited by Le Hung Son, W. Tutschke, C.C. Yang, Kluwer
Academic Publishers 2003, pp.45-83.

\medbreak\noindent[Siu 2005]  Y.-T. Siu, Multiplier ideal sheaves in
complex and algebraic geometry, {\it Science in China,
Ser.A:Math.}{\bf 48} (2005), 1-31 (arXiv:math.AG/0504259).

\medbreak\noindent[Siu 2006] Y.-T. Siu, A general non-vanishing
theorem and an analytic proof of the finite generation of the
canonical ring, arXiv:math/0610740.

\medbreak\noindent[Siu 2007] Y.-T. Siu, Additional explanatory notes
on the analytic proof of the finite generation of the canonical
ring, arXiv:0704.1940.

\medbreak\noindent[Skoda 1972] H. Skoda, Application des techniques
$L^2$ \`a la th\'eorie des id\'eaux d'une alg\`ebre de fonctions
holomorphes avec poids. {\it Ann. Sci. \'Ecole Norm. Sup.} {\bf 5}
(1972), 545-579.

\medbreak\noindent[Viehweg 1982] E. Viehweg, Vanishing theorems.
{\it J. Reine Angew. Math.} {\bf 335} (1982), 1-8.

\bigbreak\noindent{\it Author's mailing address}: Department of
Mathematics, Harvard University, Cambridge, MA 02138, U.S.A.

\medbreak\noindent {\it Author's e-mail address}:
siu@math.harvard.edu

\end{document}